\documentclass[11pt,english,uppersorbian,american]{article}
\usepackage{bm}%
\usepackage[colorlinks=true,linkcolor=black,citecolor=black,urlcolor=black]{hyperref}%

\usepackage[T1]{fontenc}
\usepackage[utf8]{inputenc}
\usepackage{algorithm}
\usepackage[font={small,it}]{caption}
\usepackage{authblk}
\usepackage{float}
\usepackage{mathtools}
\usepackage{amsmath}
\usepackage{amsthm}
\usepackage{amssymb}
\usepackage{graphicx}
\usepackage{algpseudocode}
\usepackage[round]{natbib}
\usepackage{appendix}

\makeatother
\usepackage{subfig}

\floatstyle{ruled}
\newfloat{algorithm}{tbp}{loa}
\providecommand{\algorithmname}{Algorithm}
\floatname{algorithm}{\protect\algorithmname}

\theoremstyle{definition}
\newtheorem{defn}{\protect\definitionname}
\theoremstyle{plain}

\theoremstyle{remark}

\providecommand{\definitionname}{Definition}

\providecommand{\theoremname}{Theorem}

\makeatother

\providecommand{\definitionname}{Definition}
\providecommand{\remarkname}{Remark}
\providecommand{\theoremname}{Theorem}

\expandafter\ifx\csname package@font\endcsname\relax\else
\expandafter\expandafter
\expandafter\usepackage
\expandafter\expandafter
\expandafter{\csname package@font\endcsname}%
\fi
\author{Mattia Serra\thanks{serram@ethz.ch}}
\author{George Haller\thanks{Email address for correspondence: georgehaller@ethz.ch}}
\title{Forecasting Long-Lived Lagrangian Vortices from their Objective Eulerian Footprints}%
\affil{Institute for Mechanical Systems, ETH Zürich\\Leonhardstrasse 21, 8092 Zurich, Switzerland}%

\begin{document}

	\maketitle
\begin{abstract}
\noindent We derive a non-dimensional metric to quantify the expected Lagrangian
persistence of objectively defined Eulerian vortices in two-dimensional
unsteady flows. This persistence metric is the averaged deviation
of the vorticity from its spatial mean over the Eulerian vortex, normalized
by the instantaneous material leakage from the Eulerian vortex. The
metric offers a model- and frame-independent tool for uncovering the
instantaneous Eulerian signature of long-lived Lagrangian vortices.
Using satellite-derived ocean velocity data, we show that Lagrangian vortex-persistence predictions by our
metric significantly outperform those
inferred from other customary Eulerian diagnostics, such as the potential
vorticity gradient and the Okubo-Weiss criterion.
\end{abstract}

\section{Introduction}

Coherent Lagrangian vortices \citep{LCSHallerAnnRev2015} are fluid
masses enclosed by material boundaries that exhibit only moderate
deformation under advection. Such vortices play a fundamental role
in a number of transport and mixing processes. For instance, coherent
mesoscale oceanic eddies are known to carry water over long distances,
influencing global circulation and climate \citep{BealNature2011}. 

Frame-invariant methods for the precise identification of coherent
Lagrangian vortex boundaries are now available \citep{BlackHoleHaller2013,Farazmand2015,HallerLAVD2015}.
These methods, as any Lagrangian approach, are intrinsically tied
to a preselected finite time interval. Some material vortex boundaries lose their coherence immediately beyond their extraction times, while others remain coherent over much longer intervals \citep{BlackHoleHaller2013,Beron-VeraDatassh2013,Wang2016}.
It is, therefore, of interest to identify a signature of long-lived
Lagrangian vortices without an a priori knowledge of their time scale
of existence. 

The question we address in the present paper is the following: What
instantaneous Eulerian features of a coherent Lagrangian vortex make
it likely to persist over longer time intervals? This question is
relevant, for instance, in environmental forecasting and decision-making,
as well as in assessing the life stage of coherent eddies that influence
the general circulation in the ocean. Despite its importance, however, the question
of Lagrangian vortex persistence has received little attention. Broadly
used Eulerian vortex detection methods provide no direct answer, although
the motivation for these Eulerian methods is often precisely the need
to capture sustained material transport by vortices. Clearly, the
future of advected water masses in an unsteady flow cannot be precisely predicted based on just present data. 
Reasons for this include unforeseeable future interactions with other vortices, and a priori unknown external forcing on the flow. The most one can hope for, therefore, is to forecast Lagrangian eddy persistence, with high enough probability, assuming that these unpredictable effects do not arise.

To this end, we propose here a non-dimensional metric to assess the persistence of Eulerian vortices encircled by elliptic Objective Eulerian
Coherent Structures (OECSs), as defined by \citet{SerraHaller2015}. Such OECSs are closed curves with no
short-term unevenness in their material deformation rates (zero
short-term filamentation). The objectivity of
OECSs ensures the frame-invariance of the transport estimates they
provide, while the non-dimensionality of the persistence metric introduced here will allow for a comparison of
coexisting vortices of various sizes and times scales.

Our persistence metric is the ratio of the rotational coherence strength
of an elliptic OECS to its material leakage. Eulerian vortices with
high rotation rates and low material leakage will
have high persistence metric values and will be seen to delineate regions of
sustained material coherence. As a side result, we also derive an
explicit formula for the material flux through an elliptic OECS. This
technical result is generally applicable to estimating the deformation
of limit cycles in a two-dimensional vector field under a change in
the system parameters. 

We illustrate our results on an unsteady satellite altimetry-based
velocity field of the South Atlantic Ocean. Remarkably, we find that
elliptic OECSs with high values of the persistence metric capture,
with high probability, the signature of long-lived Lagrangian vortices.
At the same time, the predictive power of customary Eulerian diagnostics,
such as the Okubo-Weiss ($OW$) criterion \citep{Okubo1970,Weiss1991},
the potential vorticity ($PV$) and the potential vorticity gradient
($\nabla PV$) \citep{Griffa2007}, turns out to be substantially lower, showing correlations below 0.5 with the actual
lifetime of Lagrangian eddies.

\section{Set-up and notation}

We consider an unsteady velocity field $v(x,t)$ defined on a spatial
domain $U\subset\mathbb{R}^{2}$ over a finite time interval $[t_{0},t_{1}]$.
We recall the velocity gradient decomposition 
\begin{equation}
\nabla v(x,t)=S(x,t)+W(x,t),\label{VelGradDec}
\end{equation}
where $S=\tfrac{1}{2}(\nabla v+\nabla v^{\top})$ and $W=\tfrac{1}{2}(\nabla v-\nabla v^{\top})$
are the rate-of-strain tensor and the spin tensor, respectively.

The spin tensor $W$ is skew-symmetric while $S$ is symmetric, with
its eigenvalues $s_{i}(x)$ and eigenvectors $e_{i}(x)$ satisfying
\[
e_{2}Se_{i}=s_{i}e_{i},\ \ \lvert e_{i}\rvert=1,\ \ i=1,2;\ \ s_{1}\leq s_{2},\ \ e_{2}=Re_{1}=\begin{bmatrix}0 & -1\\
1 & 0
\end{bmatrix}e_{1}.
\]
Fluid particle trajectories generated by $v(x,t)$ are solutions of
the differential equation $\dot{{x}}=v(x,t)$, defining the flow map

\[
F_{t_{0}}^{t}(x_{0})=x(t;t_{0},x_{0}),\ \ \ \ x_{0}\in U,\ \ \ \ t\in[t_{0,}t_{1}],
\]
which maps initial particle positions $x_{0}$ at time $t_{0}$ to
their time-$t$ positions, $x(t;t_{0},x_{0})$.

\section{Vortices as elliptic OECSs}

A typical set of fluid particles is subject to significant stretching
under advection in an unsteady flow. Even in the limit of zero advection
time, fluid elements generally experience considerable stretching
rates. One may look for the Eulerian signatures of coherent material
vortices as exceptional sets of fluid trajectories that defy this
general trend. Specifically, \citet{SerraHaller2015}
seek boundaries of Eulerian coherent vortices as closed instantaneous curves
across which the averaged material stretching rate shows no leading-order
variability.

\begin{figure*}
	\centering
	\subfloat[\label{MatbeltsInst}]{\includegraphics[width=0.47\textwidth]{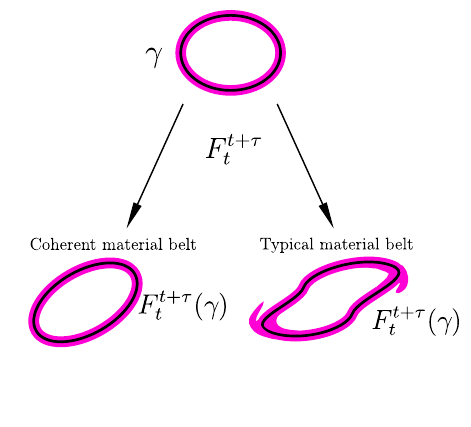}}\quad
	\subfloat[\label{fig:Globxi2t6}]{\includegraphics[width=.47\columnwidth]{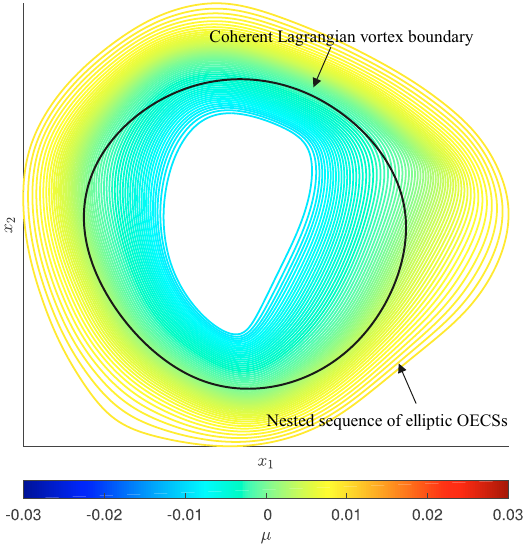}}

\caption{(a) A closed material curve $\gamma$ (black) at time $t$ is advected
by the flow into its later position $F_{t}^{t+\tau}(\gamma)$, with
$\tau\approx0$. The advected curve remains coherent if an initially
uniform material belt (magenta) around it shows no leading-order variations
in stretching rate. (b) Nested family of elliptic OECSs in a flow
example, analyzed in more detail in section \ref{PredonOcean}, for
different values of $\mu$ (in color). The elliptic OECS family fills
a region that also turns out to contain a persistent Lagrangian vortex
in this example \citep{BlackHoleHaller2013}. }
\end{figure*}

Mathematically, this is equivalent to seeking closed curves $\gamma$
whose $\mathcal{O}(\epsilon)$ perturbations show no $\mathcal{O}(\epsilon)$
variability in the averaged strain-rate functional $\dot{Q}_{t}(\gamma)$,
defined as 
\[
\dot{Q}_{t}(\gamma)=\dfrac{1}{\sigma}\oint_{\gamma}\dfrac{\langle x'(s),S(x(s),t)x'(s)\rangle}{{\langle x'(s),x'(s)\rangle}}ds.
\]
Here $x(s)$, $s\in[0,\sigma]$, denotes the arclength parametrization
of $\gamma$ at time $t$, and $x'(s)$ denotes its local tangent
vector. Stationary curves of $\dot{Q}_{t}(\gamma)$ are cores of exceptional
material belts showing perfect short-term coherence (Fig. \ref{MatbeltsInst}).
\citet{SerraHaller2015} show that closed stationary curves of
$\dot{Q}_{t}(\gamma)$ are precisely the closed null-geodesics of a suitably
defined Lorentzian metric. Along these curves, the tangential stretching
rate $\mu$ is constant. 

The closed stationary curves of $\dot{Q}_{t}(\gamma)$ turn out to
be computable as limit cycles of the direction field family 
\begin{equation}
x'=\chi_{\mu}^{\pm}(x),\ \ \chi_{\mu}^{\pm}(x)=\sqrt{\dfrac{s_{2}(x)-\mu}{s_{2}(x)-s_{1}(x)}}e_{1}(x)\ \pm\ \sqrt{\dfrac{\mu-s_{1}(x)}{s_{2}(x)-s_{1}(x)}}e_{2}(x),\label{etamufield}
\end{equation}
within the domain $U_{\mu}\subset U$ defined as 
\[
U_{\mu}=\{x\in U\lvert\ s_{2}-s_{1}\neq0,\ s_{1}\le\mu\le s_{2}\}.
\]
The direction field family (\ref{etamufield}) depends on the choice
of the sign parameter $\pm$, as well as on the parameter $\mu\in\mathbb{R}$.
We define elliptic OECSs as limit cycles of (\ref{etamufield}) for
each value of the parameter $\mu\approx0$. The $\mu=0$ member of
this one-parameter family of nested curves represents a perfect instantaneously
coherent vortex boundary \citep{SerraHaller2015}. Such a closed curve
is highly atypical, exhibiting no instantaneous stretching rate.

Members of limit cycles families of $\chi_{\mu}^{\pm}$ cannot intersect.
Each limit cycle either grows or shrinks under changes in $\mu$,
forming a smooth annular belt of non-intersecting loops (see \citet{SerraHaller2015}
for details). This annular Eulerian belt often surrounds a persistent
Lagrangian vortex boundary, as in the example shown in Fig. \ref{fig:Globxi2t6}.

\section{Material flux through elliptic OECSs}

In this section, we derive an explicit formula for the material flux
through an elliptic OECS to quantify the degree to which the OECS
is Lagrangian. As a byproduct, we obtain an expression for the short-term
continuation of elliptic OECSs under varying time. 

Let $\gamma(t)$ be a time-varying, closed curve family parametrized
by a function $x(s,t)$. The pointwise instantaneous material flux
density through $\gamma(t)$ is then given by

\begin{equation}
\begin{aligned}\varphi(x(s,t),t)= & \big\langle v(x(s,t),t)-\tfrac{d}{dt}x(s,t),\ n(x(s,t),t)\big\rangle\\
= & \big\langle v(x(s,t),t),\ n(x(s,t),t)\big\rangle-\big[\tfrac{d}{dt}x(s,t)\big]^{\perp},
\end{aligned}
\label{PointwiseInstFlux}
\end{equation}
i.e., by the curve-normal projection $\langle\cdot,n(x(s,t),t)\rangle$
of the Lagrangian velocity $v(x(s,t),t)$ of a trajectory relative
to the velocity of $\gamma(t)$.

In our context, $x(s,t)$ represents a limit cycle of the ODE (\ref{etamufield}),
thus we have $n(x(s,t),t)=[\chi_{\mu}^{\pm}(x(s,t),t)]^{\perp}=R\chi_{\mu}^{\pm}(x(s,t),t)$.
In \ref{AppFlux}, we derive and solve an ODE for the unknown
term $\big[\tfrac{d}{dt}x(s,t)\big]^{\perp}$ in (\ref{PointwiseInstFlux}),
obtaining the final formula 
\begin{equation}
\big[\tfrac{d}{dt}x(s,t)\big]^{\perp}=\Phi_{0}^{s}(t)\big[\tfrac{d}{dt}x(0,t)\big]^{\perp}+\Pi(s,t),\ \ \Pi(s,t):=\Phi_{0}^{s}(t)\int_{0}^{s}(\Phi_{0}^{\vartheta}(t))^{-1}\tilde{c}(\vartheta,t)d\vartheta,\label{VArOfConst}
\end{equation}
with $\Phi_{0}^{s}(t)$ denoting the matrix 
\begin{equation}
\begin{split}\Phi_{0}^{s}(t) & =\begin{bmatrix}1 & \int_{0}^{s}e^{\int_{0}^{\vartheta}\nabla\cdot\chi_{\mu}^{\pm}(x(\vartheta,t),t)d\vartheta}\kappa(x(\vartheta,t))d\vartheta\\
0 & e^{\int_{0}^{s}\nabla\cdot\chi_{\mu}^{\pm}(x(\vartheta,t),t)d\vartheta}
\end{bmatrix},\end{split}
\label{fundmatrsolinetabasesimplified}
\end{equation}
and 
\begin{equation}
\begin{aligned}\tilde{c}(s,t) & =[0,\psi(x(s,t),t)]^{\top},\\
\psi(x(s,t),t) & =\frac{-\langle\chi_{\mu}^{\pm}(s,t),\partial_{t}S(s,t)\chi_{\mu}^{\pm}(s,t)\rangle}{2\langle\chi_{\mu}^{\pm}(s,t),S(s,t)\chi_{\mu}^{\pm}(s,t)^{\perp}\rangle},\\
\kappa(x(s,t)) & =\langle\nabla\chi_{\mu}^{\pm}(s,t)\chi_{\mu}^{\pm}(s,t),R\chi_{\mu}^{\pm}(s,t)\rangle.
\end{aligned}
\label{Partialt0etavFt}
\end{equation}
Note that $\kappa$ represents the pointwise curvature along the elliptic
OECS with respect to the normal vector defined as $[\chi_{\mu}^{\pm}]^{\perp}=R\chi_{\mu}^{\pm}$.

In \ref{AppFlux}, we also derive the following equation
for the correct initial condition of  $\big[\tfrac{d}{dt}x(s,t)\big]^{\perp}$:
\begin{equation}
\big[\tfrac{d}{dt}x(0,t)\big]^{\perp}=\dfrac{\langle\Pi(\sigma,t),d\rangle}{1-\rho_{2}(t)},\ \ \ \rho_{2}(t)=e^{\int_{0}^{\sigma}\nabla\cdot\chi_{\mu}^{\pm}(x(\vartheta,t),t)d\vartheta},\ \ \ d:=\left[\begin{array}{c}
0\\
1
\end{array}\right].\label{z0final}
\end{equation}

\noindent This initial condition represents the ratio between the
magnitude of the perturbation needed to destroy the limit cycle due
to the unsteadiness of the flow and the strength of the hyperbolicity
of the limit cycle. For steady flows, we have $\Pi^{\perp}(s,t)=0$
since $\tilde{c}(s,t)=0$. In that case, the robustness of the limit
cycle is determined by $\rho_{2}(t)\equiv\rho_{2}=const.$, without any
time dependence.

For the computation of $\partial_{t}S$ in (\ref{Partialt0etavFt}),
we employ a backward finite-difference scheme (see \ref{SecAutomScheme}
for details). Once $\big[\tfrac{d}{dt}x(0,t)\big]^{\perp}$ is known,
we evaluate the pointwise flux density introduced in (\ref{PointwiseInstFlux}).
For a counterclockwise parametrization of $\gamma(t)$, and for our
definition of $[\chi_{\mu}^{\pm}]^{\perp}$, positive values of $\varphi(x(s,t),t)$
represents inward material flux.

\begin{figure}
\centering \includegraphics[width=0.8\textwidth]{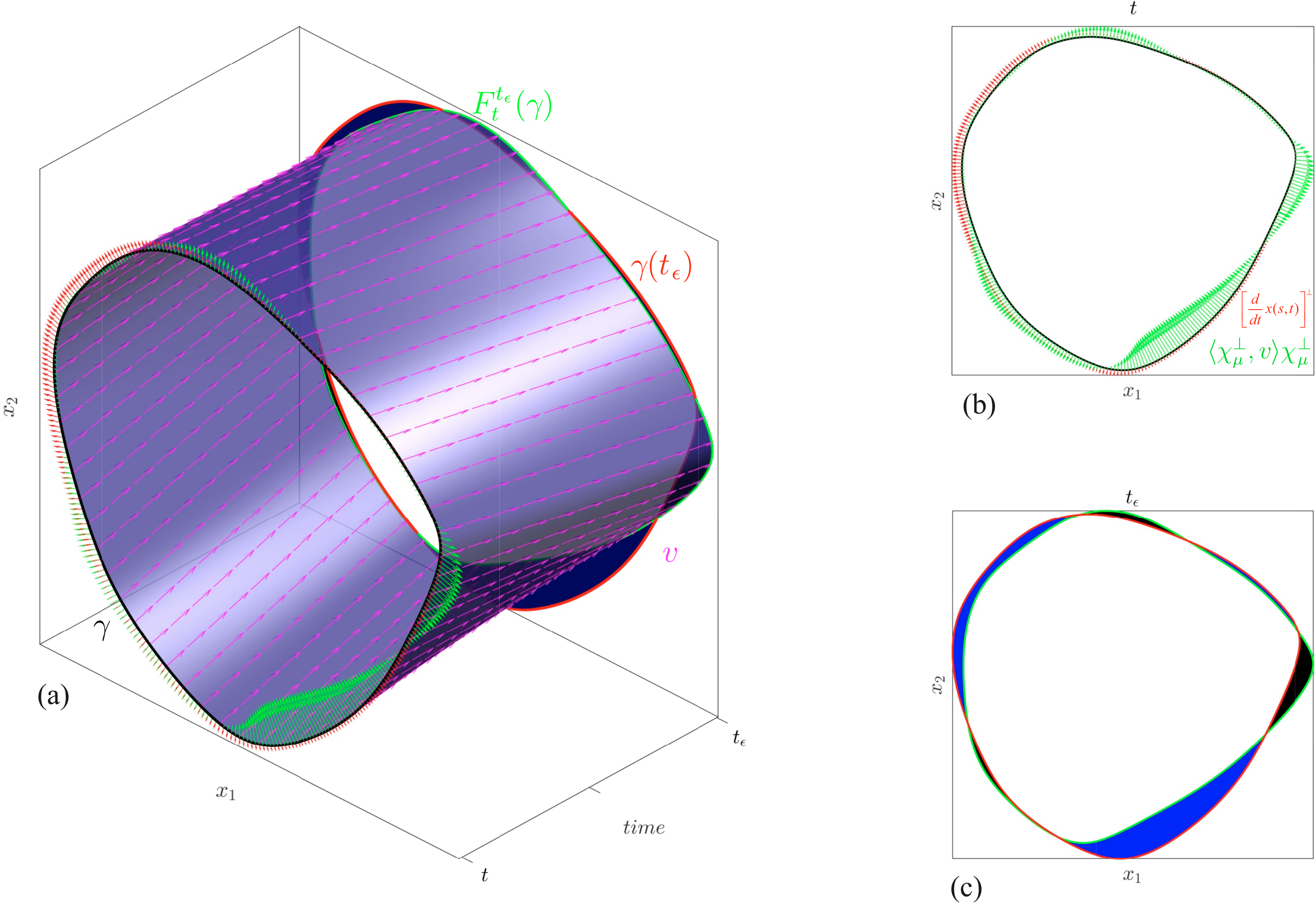}
\caption{(a) Initial elliptic OECS $\gamma$ (black) and its advected image
under the flow in the extended phase space over $[t,\ t_{\epsilon}]$,
where $t_{\epsilon}=t+\epsilon\Delta t$. At time $t$, the flow velocity
perpendicular to the curve and the corresponding elliptic OECS velocity
are reported by the green and red arrows, respectively. At time $t_{\epsilon}$,
the advected image, $F_{t}^{t_{\epsilon}}(\gamma)$, is shown in green
while the recomputed elliptic OECS $\gamma(t_{\epsilon})$ in red.
(b) Slice of (a) at time $t$. (c) Slice of (a) at time $t_{\epsilon}$.
The blue and black areas represent the actual inward and outward material
flux across $\gamma$ over $[t,t_{\epsilon}]$, respectively.}
\label{Flux3dfig} 
\end{figure}

Figure \ref{Flux3dfig}a illustrates an elliptic OECS $\gamma(t)$
(black) at time $t$, for a fixed value of $\mu$, with its advected
image over the time window $[t,\ t_{\epsilon}]$ in the extended phase
space of position and time. Figures \ref{Flux3dfig}a-b show the initial
and final time slices of Fig. \ref{Flux3dfig}a. The materially advected
image of $\gamma$ at time $t_{\epsilon}$, $F_{t}^{t_{\epsilon}}(\gamma)$,
is shown in green while the elliptic OECS $\gamma(t_{\epsilon})$
computed at time $t_{\epsilon}$, is shown in red. Figures \ref{Flux3dfig}
a-b show the instantaneous pointwise material flux density through
$\gamma(t)$, given by the difference between the flow velocity normal
to the curve (green arrows) and the corresponding continuation velocity
(red arrows). Given formula \eqref{PointwiseInstFlux}, the total
instantaneous material flux across $\gamma(t)$ is 
\[
\varphi_{\gamma(t)}=\oint_{\gamma(t)}\varphi(x(s,t),t)ds,
\]
with 

\[
\begin{cases}
\begin{aligned}\varphi(x(s,t),t)= & \langle v(x(s,t),t),[\chi_{\mu}^{\pm}]^{\perp}(x(s,t),t)\rangle-\langle\Phi_{0}^{s}(t)\big[\tfrac{d}{dt}x(0,t)\big]^{\perp}+\Pi(s,t),d\rangle,\\
\big[\tfrac{d}{dt}x(0,t)\big]^{\perp}= & \dfrac{\langle\Pi(\sigma,t),d\rangle}{1-e^{\int_{0}^{\sigma}\nabla\cdot\chi_{\mu}^{\pm}(x(\vartheta,t),t)d\vartheta}}.
\end{aligned}
\end{cases}
\]
The instantaneous total material flux $\varphi_{\gamma(t)}$, multiplied
by $\epsilon\Delta t$, approximates the actual material flux given
by the inward (blue) area minus the outward (black) area shown in
Fig. \ref{Flux3dfig}c.

\section{Persistence metric for elliptic OECSs}

We now propose a non-dimensional objective metric that classifies
elliptic OECSs based on their expected persistence in time. We first
define the two ingredients needed for this metric: the \textit{rotational
coherence} and the \textit{relative material leakage.} 
\begin{defn}
The rotational coherence of an elliptic OECS $\gamma(t)$ is
\begin{equation}
\omega_{\gamma}(t):=\dfrac{\vert\int_{A_{\gamma(t)}}[\omega(x,t)-\overline{\omega}(t)]dA\vert}{A_{\gamma(t)}},\label{ObjMeanVort}
\end{equation}
where $\omega(x,t)$ denotes the vorticity, $A_{\gamma(t)}$ is the
area enclosed by $\gamma(t)$, and 
\[
\overline{\omega}(t)=\dfrac{\int_{A_{\partial U}}\omega(x,t)dA}{A_{\partial U}}
\]
is the mean spatial vorticity over the domain $U$ with boundary
$\partial U$. 
\end{defn}
The rotational coherence $\omega_{\gamma}$ represents the normed
mean vorticity deviation within $\gamma(t)$, inspired by related
quantities defined in \citet{HallerLAVD2015}. Specifically, the rotational
coherence measures the strength of a vortical structure arising from
its rotational speed. The classic measure of vortex strength, also
called circulation \citep{Batchelor2000}, relies solely on the vorticity
$\omega(x,t)$, and is therefore frame-dependent. The rotational coherence
$\omega_{\gamma}$, instead, involves the vorticity deviation, which
is frame-independent (\ref{Objectivity}). Elliptic OECSs with high rotational coherence are shielded by locally high levels of shear, and hence are expected to persist in time.

\begin{defn}
The relative material leakage of an elliptic OECSs $\gamma(t)$ is 
\begin{equation}
\Gamma_{\gamma}(t):=\dfrac{\oint_{\gamma(t)}\lvert\varphi(x(s,t),t)\rvert ds}{A_{\gamma(t)}}.\label{AbsoluteFlux}
\end{equation}

\end{defn}
The relative material leakage measures the rate of material area leaking
out of $\gamma(t)$ due to its non-Lagrangian evolution, divided by
the initial area of $\gamma(t)$. A $\gamma(t)$ with low $\Gamma_{\gamma}(t)$
identifies an exceptional curve that exhibits low inhomogeneity in
its stretching rates both in its initial position and in its short-term
advected position. The absolute value in \eqref{AbsoluteFlux} prevents
the cancellation of opposite-sign material flux contributions. Note
that both $\omega_{\gamma}$ and $\Gamma_{\gamma}$ have the dimension
$[time^{-1}]$. 

We expect elliptic OECSs with high rotational coherence and low material leakage to be the best candidate locations for Lagrangian
vortices. To this end, we define the\textit{ persistence metric} of
an elliptic OECS as the following objective, non-dimensional quantity: 
\begin{defn}
\label{DefVortexPersistence} The persistence metric of elliptic OECS
$\gamma(t)$ is
\begin{equation}
\Theta_{\gamma}(t):=\dfrac{\textit{ rotational coherence}}{\textit{relative material leakage }}=\dfrac{\omega_{\gamma}(t)}{\Gamma_{\gamma}(t)}=\dfrac{\vert\int_{A_{\gamma(t)}}[\omega(x,t)-\overline{\omega}(t)]dA\vert}{\oint_{\gamma(t)}\lvert\varphi(x(s,t),t)\rvert ds}.\label{NondimMetric}
\end{equation}
\end{defn}
The non-dimensional nature of $\Theta_{\gamma}$ is immediate from
equation (\ref{NondimMetric}), while its frame-invariance follows
from the objectivity of the scalar quantities involved in its definition
(cf. \ref{Objectivity}). The non-dimensionality of $\Theta_{\gamma}$
allows us to characterize the persistence of vortices regardless of
the their spatial and temporal scales, which are often abundant and
unknown. The objectivity of $\Theta_{\gamma}$ ensures a persistence assessment independent of the frame of reference.

In case of zero relative material leakage, we have $\Theta_{\gamma}=\infty$,
as indeed desired for a perfectly material elliptic OECSs. In this
rare case, $F_{t}^{t_{\epsilon}}(\gamma(t))=\gamma(t_{\epsilon})$
and hence the green and the red curves in Fig. \ref{Flux3dfig}c
coincide. 

In \ref{SecAutomScheme}, we summarize the numerical algorithms
for the identification of likely long-lived Lagrangian vortices
from their objective Eulerian features. Specifically, Algorithm 1 summarizes
the computation of elliptic OECSs and Algorithm 2 describes the computation
of the corresponding persistence metric $\Theta_{\gamma}$.

\section{Example: Forecasting persistent Lagrangian vortices in satellite-derived
ocean velocity data\label{PredonOcean}}

We apply our OECS-based vortex-coherence forecasting scheme to a two-dimensional
unsteady ocean dataset obtained from AVISO satellite altimetry measurements
(\url{http://www.aviso.oceanobs.com}). The domain of interest is
the Agulhas leakage in the Southern Atlantic Ocean bounded by longitudes
$[17^{\circ}W,7^{\circ}E]$ and latitudes $[38^{\circ}S,22^{\circ}S]$.
The Agulhas Current is a narrow western boundary current of the southwest
Indian Ocean, whose interaction with the strong Antarctic Circumpolar
Current gives rise to Agulhas rings, the largest mesoscale eddies
in the ocean.

Agulhas rings are considered important in the global circulation due
to the large amount of water they carry over considerable distances
\citep{BealNature2011}. For comparison with earlier Lagrangian
analysis \citep{BlackHoleHaller2013}, we consider the same initial
time $t=24\ November\ 2006$ and a similar but slightly larger spatial
domain. For more detail on the dataset and the numerical method, see \ref{SecAutomScheme} and \ref{SSHdataset}. 

As mentioned earlier, the $OW$ parameter
\[
OW(x,t)=s_{2}^{2}(x,t)-\omega^{2}(x,t),
\]
is a frequently used indicator of instantaneous ellipticity in unsteady
fluid flows \citep{Okubo1970,Weiss1991}. Spatial domains
with $OW(x,t)<0$ (rotation prevailing over strain) are generally
considered vortical. The $OW$ parameter is not objective (the vorticity
term will change under rotations), and hence no objective threshold
level can be defined for this scalar field to identify vortices unambiguously.
This ambiguity significantly impacts the overall number and geometry
of the vortical structures inferred from the $OW$ parameter (\ref{OWThresholding}). Among other applications, $OW$ has been used
to study eddies in the Gulf of Alaska \citep{OWls_theshld_Henson2008},
in the Mediterranean Sea \citep{OWlsthrsh_Isern-Fontanet2006,Isern-Fontanet2004,Isern-Fontanet2003},
in the Tasman Sea \citep{Waugh2006}, and in the global ocean \citep{Chelton2007}. 

\begin{figure}
\centering \includegraphics[width=1\textwidth]{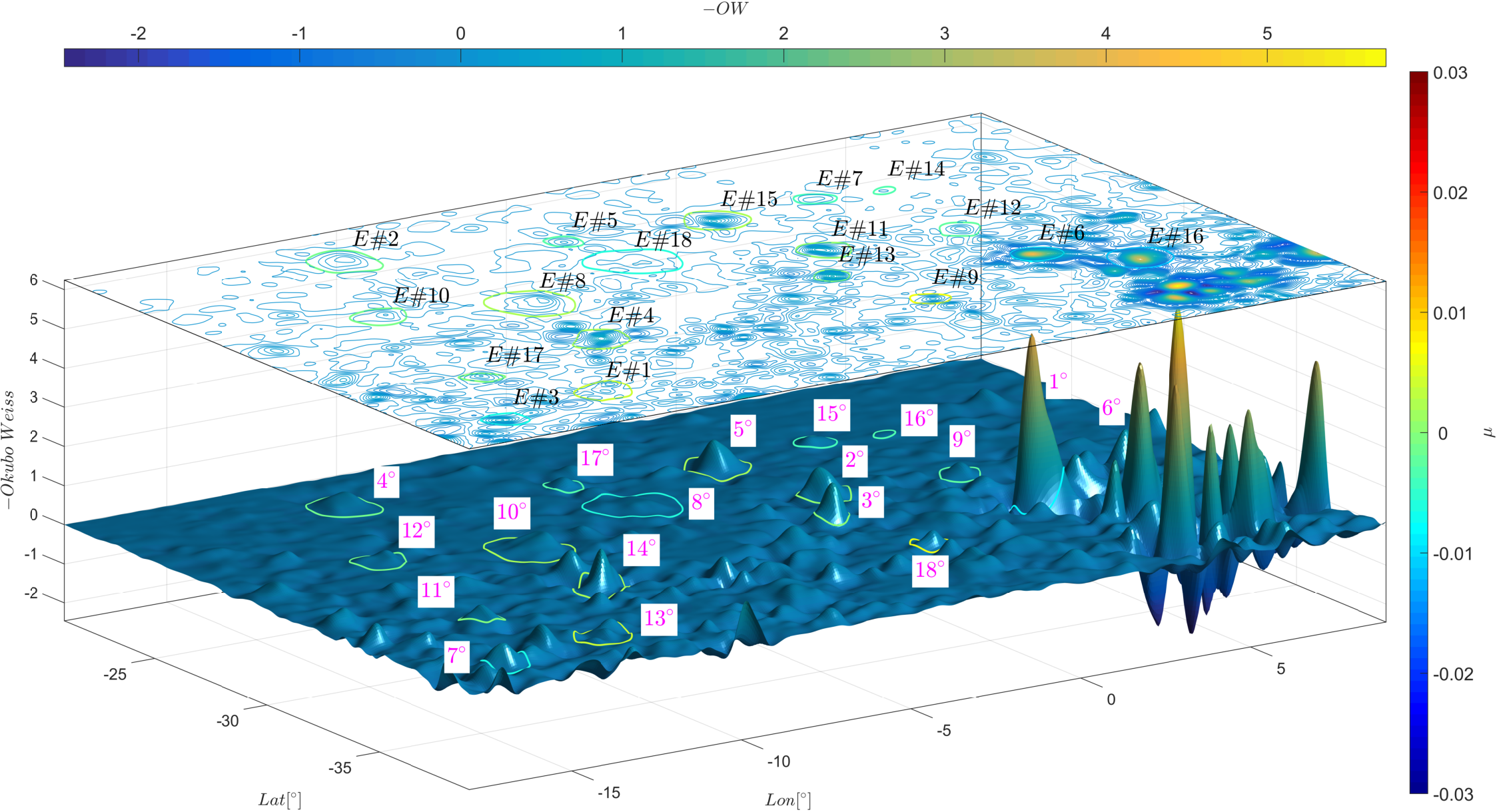}
\caption{Elliptic OECSs with the highest persistence metric $\Theta_{\gamma}$
on a surface representing the negative $OW$ parameter (horizontal
colorbar or z-axis). The color of Elliptic OECSs represents the corresponding
stretching-rate value $\mu$ (right colorbar). Black numbers identify
different vortical regions detected by elliptic OECSs. In magenta,
the classification of the most persistent vortical regions in decreasing
order of $\Theta_{\gamma}$.}
\label{BestISGonOW3dsurf} 
\end{figure}

In Fig. \ref{BestISGonOW3dsurf}, we show elliptic OECSs with the highest
persistence metric $\Theta_{\gamma}$ for each vortical region, on
a surface representing the negative $OW$ parameter. The plane of
the figure also shows the level curves of the $OW$ parameter. The
black numbers in Fig. \ref{BestISGonOW3dsurf} label the different
vortical structures, while the magenta numbers classify them in decreasing
order of $\Theta_{\gamma}$. We find elliptic OECSs in locations of
the flow where the $OW$ parameter is close to zero and hence signals
no vortices (see, e.g., $E\#7$, $E\#8$, $E\#18$). In contrast, close
to the tip of Africa, $OW$ signals several strong vortical regions,
even though we only detect two belts of elliptic OECSs ($E\#6$, $E\#16$).

To assess these discrepancies between the $OW$ parameter and our
persistence metric $\Theta_{\gamma}$, we compare the coherence strength
suggested by $\Theta_{\gamma}$ to the actual lifetime of Lagrangian
vortices computed over a time window of four months with initial time
$t=24\ November\ 2006$. We compute the Lagrangian lifetime of elliptic
OECS's as the maximum integration time for which coherent (i.e., non-filamenting)
Lagrangian vortices in the sense of \citet{BlackHoleHaller2013} exist
nearby. To do so, we use the following discrete set of integration
times: $[7,15,30,60,90,120]$ days. The Lagrangian lifetime of an
elliptic OECSs is then the largest integration time from this sequence that still gives a nearby coherent Lagrangian eddy. We consider a Lagrangian eddy to be in the vicinity of an elliptic
OECSs if it is contained within a circle of radius $3^{\circ}$ ($\sim$1.5
times the radius of a mesoscale eddy) centered at the elliptic OECSs.
\begin{figure}
\centering \includegraphics[width=0.9\textwidth]{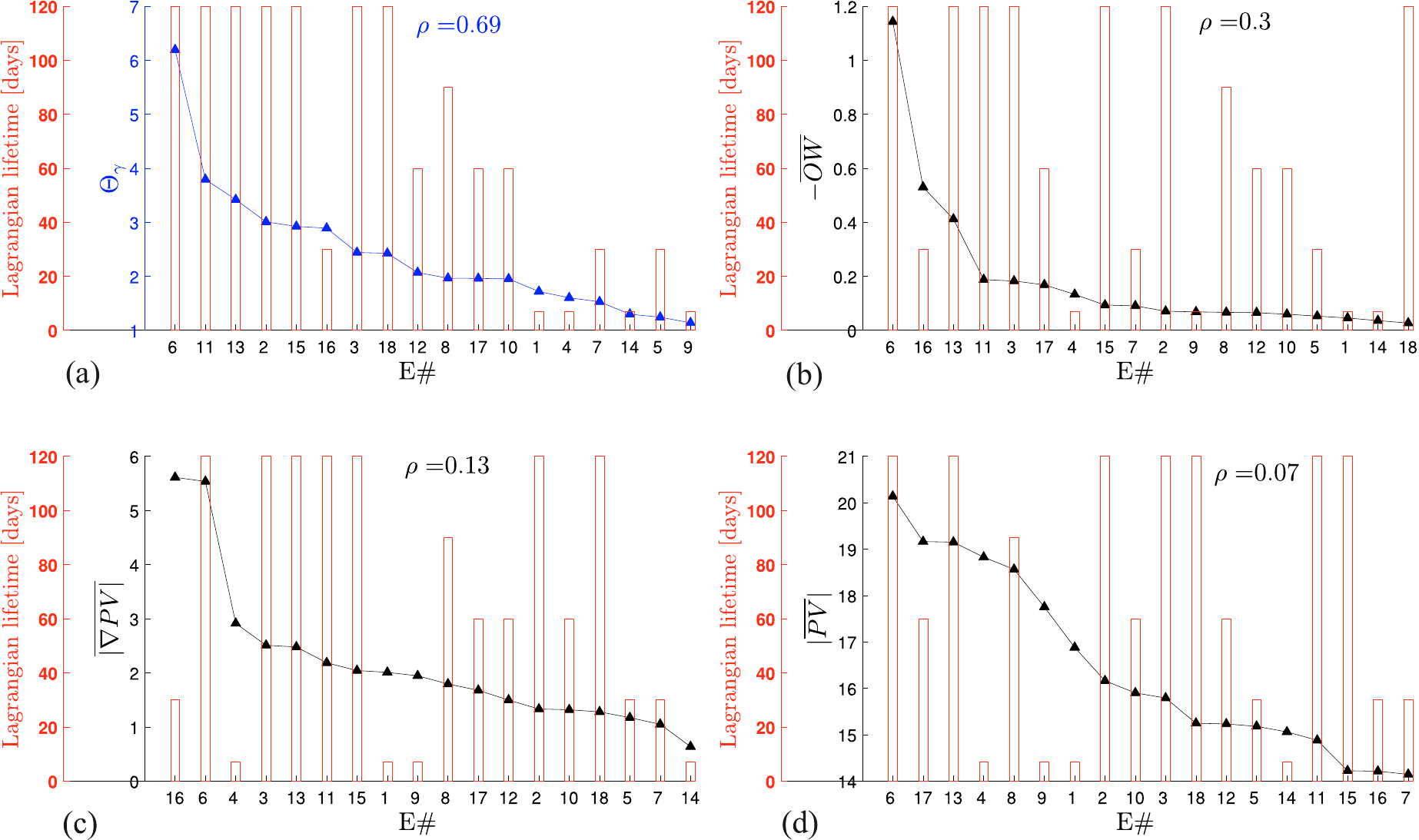}
\caption{(a) Values of the persistence metric $\Theta_{\gamma}$ (blue) for
the different vortical regions identified by elliptic OECSs compared
with their Lagrangian lifetime (red). (b-d) Spatial average of $-OW$,
$\vert\nabla PV\vert$ and $\vert PV\vert$ within elliptic OECSs
compared with their Lagrangian lifetime (red). The parameter $\rho$
indicates the correlation coefficient between the instantaneous prediction
given by each metric and the actual Lagrangian lifetime of the underlying
vortical region.}
\label{MetricGLMOW} 
\end{figure}

Figure \ref{MetricGLMOW}a shows the $\Theta_{\gamma}$ values (blue)
associated with each vortical region $(E\#i)$ in descending order
of $\Theta_{\gamma}$. Figures \ref{MetricGLMOW}b-d, in contrast,
show alternative instantaneous metrics, such as the average of $-OW$,
$\vert\nabla PV\vert$ and $\vert PV\vert$, respectively, within
the elliptic OECSs shown in Fig. \ref{BestISGonOW3dsurf}. For this
dataset, we compute $PV$ as in \citet{Early2011}. The actual Lagrangian
lifetime of the underlying vortical regions is shown in red in all
the plots, along with its correlations with the different instantaneous
metrics.

The instantaneous persistence metric $\Theta_{\gamma}$ shows a distinct
correlation ($\rho\approx0.7$) with the lifetime of long-lived Lagrangian eddies in our
study domain. This includes eddies $\#6,\#11,\#13,\#15,\#18,\#8,\#23$
and $\#3$, previously identified as exceptionally coherent Lagrangian
eddy regions in \citet{BlackHoleHaller2013} and \citet{AutoMdetectFlorian}.
Figure \ref{MetricGLMOW}a shows that out of the ten elliptic OECSs
with the highest $\Theta_{\gamma}$ values, eight are long-lived Lagrangian
vortices. 

At the same time, Figs. \ref{MetricGLMOW}b-d, reveal weak predictive
power for other instantaneous Eulerian diagnostics, each of which
has significantly lower correlation with the Lagrangian lifetime of
eddies. Figures\ref{MetricGLMOW}b-c show that long-lived mesoscale
eddies, such as $E\#2$ and $E\#18$, have surprisingly weak signatures
in the $OW$ and $\vert\nabla PV\vert$ fields, while the vortex $\#16$,
which has a relatively low Lagrangian lifetime, has the strongest
signature in these two fields. The correlation coefficients of these
diagnostics would be even lower if the candidate vortical regions
were identified from the usual ad hoc threshold values for these methods,
instead of elliptic OECSs. Indeed, note that within all the vortical
regions in the southeast of the domain signaled by $OW$, only one
($E\#6$) predicts correctly a long-lived mesoscale eddy.

Regions of high $PV$ gradient are also frequently used as indicators
of instantaneous ellipticity in unsteady fluid flows. Accordingly,
in \ref{OWThresholding}, we plot the Elliptic OECSs of Fig.
\ref{BestISGonOW3dsurf} again over the $\vert\nabla PV\vert$ scalar
field. Similarly to the the $OW$-criterion, the $\vert\nabla PV\vert$
diagnostic highlights regions where no long-lived Lagrangian eddies
are present, while it misses regions where such eddies are known to
be present. 

One may alternatively compute the Lagrangian lifetime of eddies from
other objective elliptic LCS detection methods, such as the Polar
Rotation Angle (PRA) defined by \citet{Farazmand2015} or the Lagrangian-averaged
vorticity deviation (LAVD) introduced by \citet{HallerLAVD2015}. The results (not
shown here) obtained in this fashion are close to those in Fig. \ref{MetricGLMOW}.

\section{Conclusions}

We have introduced a frame-invariant, non-dimensional metric to assess
the ability of elliptic objective Eulerian coherent structures (OECS)
to identify vortical regions with sustained material coherence. Our
metric $\Theta_{\gamma}$ is the ratio between a rotational coherence
measure of the vortex and the material leakage out of the vortex. 

We have tested the $\Theta_{\gamma}$ metric on satellite-derived
ocean velocity data, where we found that Elliptic OECSs with high
$\Theta_{\gamma}$ values tend to forecast the exceptionally coherent
Lagrangian vortices found in \citet{BlackHoleHaller2013} with high
probability. To our knowledge, this is the first Eulerian eddy census
method that is shown to display a clear correlation with the actual
lifetime of nearby Lagrangian vortices. In contrast, we have found
other available Eulerian vortex diagnostics to show a distinct lack
of correlation with long-term Lagrangian coherence. This is perhaps
unsurprising because none of them is non-dimensional or objective,
and none of them is inferred from the infinitesimally short-time limit
of a mathematically exact Lagrangian coherence criterion. The lack
of correlation of classic Eulerian vortex diagnostics with Lagrangian
eddy lifetimes is consistent with the findings of \citet{Beron-VeraDatassh2013} and \citet{Wang2016}, who show that these diagnostics overestimate
the number of materially coherent vortices significantly.

Our proposed vortex persistence metric is purely kinematic,
and hence offers a model-independent instantaneous forecasting tool.
This tool is free from kinetic assumptions, such as conservation or
near-conservation of vorticity or potential vorticity.

Based on the results presented here, we expect our approach to be
useful in real-time transport predictions, environmental decision
making and hazard assessment. The purpose of this
study has been to demonstrate the predictive power of the proposed
persistence metric. A more detailed statistical analysis is planned for future work.

\newpage
\appendix
\gdef\thesection{Appendix \Alph{section}}

\section{Material flux through elliptic OECSs}\label{AppFlux}

Here we derive a formula for the instantaneous material flux through
an elliptic OECS $\gamma(t)$, whose arclength parametrization is
denoted by $x:s\mapsto x(s)$, with $s\in[0,\sigma]\subset\mathbb{R}$.
The closed curve $\gamma(t)$ is a limit cycle of (\ref{etamufield}),
parametrized by $s$, that depends smoothly on the time $t$. We first
observe that $\gamma(t+\epsilon\Delta t)$ persists for small $\epsilon\Delta t$.
This is guaranteed by the structural stability of limit cycles of
(\ref{etamufield}) together with the smoothness of the underlying
flow map. For small enough $\epsilon\Delta t$, therefore there exists
a nearby elliptic OECS, $\gamma(t+\epsilon\Delta t)$, that is a smooth
deformation of $\gamma(t)$.

Specifically, we can locally represent the perturbed limit cycle as
\begin{equation}
\begin{aligned}x(s,t+\epsilon\Delta t)= & x(s,t)+g(s,t;\epsilon\Delta t)\chi^{\perp}(x(s,t),t)\\
= & x(s,t)+\epsilon\Delta tg_{1}(s,t)\chi^{\perp}(s,t)+\mathcal{O}((\epsilon\Delta t)^{2}),
\end{aligned}
\label{t_0PerturbedCurve}
\end{equation}
where, $g$ and $g_{1}$ are two smooth scalar functions, and $\chi^{\perp}(x(s,t),t)$
is the local normal to the limit cycle at the point $x(s,t)\in\gamma(t)$.
(For notational simplicity we have used $\chi$ instead of $\chi_{\mu}^{\pm}$).
The period of the perturbed limit cycle is of the form $\sigma_{\epsilon}=\sigma+\epsilon\sigma_{1}+\mathcal{O}(\epsilon^{2})$,
leading to the periodicity condition 
\[
x(0,t+\epsilon\Delta t)=x(\sigma_{\epsilon},t+\epsilon\Delta t).
\]
Taylor expanding this expression with respect to $\epsilon$ and comparing
the $\mathcal{O}(\epsilon)$ terms gives 
\begin{equation}
\dot{x}(\sigma,t)=\dot{x}(0,t)-\chi(x(\sigma,t),t)\dfrac{\sigma_{1}}{\Delta t},\label{PeriodicityCondOrdeEpsilon}
\end{equation}
where the dot denotes the derivative with respect to $t$. This relation
shows that the difference between the perturbation to $\gamma(t)$
at $s=0$ and at $s=\sigma$, should be in the direction tangential
to the limit cycle $\gamma(t)$ in order to ensure its persistence
as a $\mathcal{C}^{1}$ closed curve.

In order to compute the term $\dot{x}(s,t)$ (in equation (\ref{PointwiseInstFlux})),
as well as the unknown quantities in (\ref{PeriodicityCondOrdeEpsilon}),
we write the equation of variations for the ODE (\ref{etamufield})
with respect to changes in the parameter $t$, leading to 
\begin{equation}
(\dot{x}(s,t))'=\nabla\chi(x(s,t),t)\dot{x}(s,t)+\partial_{t}\chi(x(s,t),t),\label{Eqvariwrt_t0}
\end{equation}
where the prime denotes the derivative with respect to the parameter
$s$. Equation (\ref{Eqvariwrt_t0}) is a non-autonomous linear ODE
for $\dot{x}(s,t)$. In the classic theory of dependence of solutions
on parameters, $\dot{x}(0,t)$ is generally zero since initial conditions
do not depend on the parameters. In the present case, however, the
initial condition, $x(0,t)$ does depend on $t$. This dependence
determines where the limit cycle is and how it deforms as $t$ varies.
We rewrite the ODE (\ref{Eqvariwrt_t0}) using the following shorthand
notation: 
\begin{equation}
y'(s)=A(s)y(s)+c(s),\label{Eqvariwrt_t0shorthand}
\end{equation}
where 
\begin{equation}
\begin{aligned}y(s)= & \dot{x}(s,t),\ \ \ A(s)=\nabla\chi(x(s,t),t),\ \ \ c(s)=\partial_{t}\chi(x(s,t),t),\end{aligned}
\label{shorthandNotation}
\end{equation}
with the time argument $t$ suppressed in $y,\ A$ and $c$ for brevity.

Note that $y(s)=\chi(x(s,t),t)$ is a solution to the homogeneous
part of (\ref{Eqvariwrt_t0shorthand}). As in \citet{HallerIacono2003},
we solve (\ref{Eqvariwrt_t0shorthand}) explicitly in the basis $[\chi(x(s,t),t),\chi^{\perp}(x(s,t),t)]$.
With the change of coordinates 
\begin{equation}
y(s)=T(s)z(s),\ \ T(s)=[\chi(x(s,t),t),\chi^{\perp}(x(s,t),t)],\label{CoordChange}
\end{equation}
(\ref{Eqvariwrt_t0shorthand}) can be written as 
\begin{equation}
z^{\prime}(s)=\tilde{A}(s)z(s)+\tilde{c}(s).\label{SimplifiedOdeinz}
\end{equation}
Substituting the change of coordinates (\ref{CoordChange}) into (\ref{Eqvariwrt_t0shorthand})
gives 
\begin{equation}
T'(s)z(s)+T(s)z'(s)=A(s)T(s)z(s)+c(s).\label{ChangeofCoord1}
\end{equation}
Since $T(s)\in SO(2)$, equation (\ref{ChangeofCoord1}) can be written
as 
\begin{equation}
z'(s)=[T^{\top}(s)A(s)T(s)-T^{\top}(s)T'(s)]z(s)+T^{\top}(s)c(s).\label{ChangeofCoord2}
\end{equation}
Using equations (\ref{shorthandNotation}-\ref{CoordChange}), we
can write $T^{\top}(s)A(s)T(s)$ and $T^{\top}(s)T'(s)$ as

\begin{equation}
\begin{split}T^{\top}(s)A(s)T(s)=\begin{bmatrix}\langle\chi,\nabla\chi\chi\rangle & \langle\chi,\nabla\chi\chi^{\perp}\rangle\\
\langle\chi^{\perp},\nabla\chi\chi\rangle & \langle\chi^{\perp},\nabla\chi\chi^{\perp}\rangle
\end{bmatrix},\end{split}
\label{dynamicMatrixinetabasis1}
\end{equation}

\begin{equation}
\begin{split}T^{\top}(s)T'(s)=\begin{bmatrix}\langle\chi,\nabla\chi\chi\rangle & \langle\chi,R\nabla\chi\chi\rangle\\
\langle\chi^{\perp},\nabla\chi\chi\rangle & \langle\chi^{\perp},R\nabla\chi\chi\rangle
\end{bmatrix}.\end{split}
\label{dynamicMatrixinetabasis2}
\end{equation}
Differentiating the identity $\langle\chi,\chi\rangle=1$ with respect
to $x$, we obtain the following relations 
\begin{equation}
\begin{aligned} & (\nabla\chi)^{\top}\chi=0,\  & \langle\chi,(\nabla\chi)^{\top}\chi\rangle=0,\\
 & \langle\chi,\nabla\chi\chi\rangle=0,\  & \nabla\chi\chi\perp\chi,\\
 & \langle\chi^{\perp},\nabla\chi\chi\rangle=\kappa,\  & R\nabla\chi\chi=-\kappa\chi,
\end{aligned}
\label{IdentitiesForetaflds}
\end{equation}
where $\kappa$ denotes the pointwise scalar curvature along the elliptic
OECS with respect to the normal vector defined as $\chi^{\perp}=R\chi$.
Substituting (\ref{dynamicMatrixinetabasis1}-\ref{IdentitiesForetaflds})
into (\ref{ChangeofCoord2}) leads to 
\[
\begin{split}\tilde{A}(s)=[T^{\top}(s)A(s)T(s)-T^{\top}(s)T'(s)]=\begin{bmatrix}0 & \kappa(s)\\
0 & \langle\chi^{\perp},\nabla\chi\chi^{\perp}\rangle
\end{bmatrix}.\end{split}
\]
The invariance property of the trace of a matrix under orthonormal
transformations implies that $Tr(\nabla\chi)=Tr(T^{\top}\nabla\chi T)$.
Recalling that $A=\nabla\chi$, and using equation (\ref{dynamicMatrixinetabasis1})
and equation (\ref{IdentitiesForetaflds}), we obtain 
\[
\begin{aligned}\nabla\cdot\chi= & Tr(\nabla\chi)\\
= & Tr(T^{\top}\nabla\chi T)\\
= & \langle\chi,\nabla\chi\chi\rangle+\langle\chi^{\perp},\nabla\chi\chi^{\perp}\rangle\\
= & \langle\chi^{\perp},\nabla\chi\chi^{\perp}\rangle,
\end{aligned}
\]
leading to the final form of $\tilde{A}(s)$: 
\begin{equation}
\begin{split}\tilde{A}(s)=\begin{bmatrix}0 & \kappa(x(s,t),t)\\
0 & \nabla\cdot\chi(x(s,t),t)
\end{bmatrix}.\end{split}
\label{dynamicMatrixinetabasisfinal}
\end{equation}
Now we derive a simplified expression for the forcing term of the
ODE (\ref{SimplifiedOdeinz}), i.e., for 
\begin{equation}
\begin{split}\tilde{c}(s) & =T^{\top}(s)c(s)=\begin{bmatrix}\langle\chi,\partial_{t}\chi\rangle\\
\langle\chi^{\perp},\partial_{t}\chi\rangle
\end{bmatrix}.\end{split}
\label{forcinginetabasis}
\end{equation}
To compute $\partial_{t}\chi$, we take the partial derivative of
the implicit ODE defining elliptic OECSs with respect to $t$ to obtain
\begin{equation}
\partial_{t}\langle\chi(x(s,t),t),[S(r,t)-\mu I]\chi(x(s,t),t)\rangle=0.\label{Partialt0eta}
\end{equation}
Dropping the arguments, we find equation (\ref{Partialt0eta}) equivalent
to 
\begin{equation}
\langle\chi,S\partial_{t}\chi\rangle=-\frac{\langle\chi,\partial_{t}S\chi\rangle}{2}.\label{Partialt0etav2}
\end{equation}
Since the direction field $\chi$ is normalized, we have $\partial_{t}\chi(x,t)\perp\ \chi(x,t)$,
and hence we can write 
\begin{equation}
\partial_{t}\chi(x,t)=\psi(x,t)\chi^{\perp}(x,t),\ \ \psi(x,t)\in\mathbb{R}.\label{detadt0perpetat0}
\end{equation}
Substituting (\ref{detadt0perpetat0}) into (\ref{Partialt0etav2})
leads to 
\begin{equation}
\psi(x(s,t),t)=-\frac{\langle\chi,\partial_{t}S\chi\rangle}{2\langle\chi,S\chi^{\perp}\rangle},\label{Partialt0etavF}
\end{equation}
which is always defined in the domain $U_{\mu}$, unless $\chi\equiv e_{i},i=1,2$,
in which case $\langle e_{i},Se_{i}^{\perp}\rangle=0$ .

We are interested in evaluating the instantaneous material flux through
elliptic OECSs. Along these curves, the constant instantaneous stretching
rate $\mu$ is approximately zero, and hence the $\chi_{\mu}^{\pm}$
directions are far from the $e_{i}$ directions. Specifically, for
incompressible flows, the directions $\chi_{0}^{\pm}$ exactly bisect
the $e_{i}$ directions. Therefore, equation (\ref{Partialt0etavF})
is always well-defined on elliptic OECSs. Substituting (\ref{detadt0perpetat0})
and (\ref{Partialt0etavF}) into (\ref{forcinginetabasis}) leads
to 
\begin{equation}
\begin{split}\tilde{c}(s) & =\begin{bmatrix}0\\
\psi(x(s,t),t)
\end{bmatrix},\end{split}
\label{forcinginetabasisF}
\end{equation}
as in (\ref{Partialt0etavFt}). 

Using the variation of constants formula (see e.g., \citet{ArnoldODE1973}),
we can write the solution of (\ref{SimplifiedOdeinz}) as 
\[
\begin{aligned}z(s)= & \Phi_{0}^{s}z(0)+\Phi_{0}^{s}\int_{0}^{s}(\Phi_{0}^{\vartheta})^{-1}\tilde{c}(\vartheta)d\vartheta\\
= & \Phi_{0}^{s}z(0)+\Pi(s),
\end{aligned}
\]
with $\Phi_{0}^{s}$ being the normalized fundamental matrix solution
to the homogeneous problem 
\begin{equation}
z'(s)=\tilde{A}(s)z(s).\label{Eqvariwrt_t0shorthandHom}
\end{equation}
By direct integration of (\ref{Eqvariwrt_t0shorthandHom}) we obtain
\[
\begin{split}\Phi_{0}^{s} & =\begin{bmatrix}1 & \int_{0}^{s}e^{\int_{0}^{\vartheta}\nabla\cdot\chi_{\mu}^{\pm}(x(\vartheta,t),t)d\vartheta}\kappa(x(\vartheta,t))d\vartheta\\
0 & e^{\int_{0}^{s}\nabla\cdot\chi_{\mu}^{\pm}(x(\vartheta,t),t)d\vartheta}
\end{bmatrix},\end{split}
\]
as in (\ref{fundmatrsolinetabasesimplified}).

Once this fundamental matrix solution is computed, the only missing
quantity in (\ref{VArOfConst}) is the initial condition $z(0)$.
To obtain that, we rewrite (\ref{PeriodicityCondOrdeEpsilon}) in
the $z$ coordinates. This, together with (\ref{VArOfConst}), leads
to the system
\begin{equation}
\begin{cases}
z(\sigma)=z(0)-d\dfrac{\sigma_{1}}{\Delta t},\ \ \ d:=[0,1]^{\top}\\
z(\sigma)=\Phi_{0}^{\sigma}z(0)+\Pi(\sigma).
\end{cases}\label{systy0}
\end{equation}
Although this system of equations is undetermined ($z(\sigma),\ z(0)$
and $\sigma_{1}$ are unknown), it is sufficient to determine the
component of $z(0)$ along the $\chi^{\perp}$ direction, $z^{\perp}(0)$.
Substituting (\ref{VArOfConst}) and (\ref{fundmatrsolinetabasesimplified})
into (\ref{systy0}), we obtain 
\[
\begin{split}\begin{bmatrix}0 & \int_{0}^{\sigma}e^{\int_{0}^{y}\nabla\cdot\chi_{\mu}^{\pm}(x(\vartheta,t),t)d\vartheta}\kappa(x(y,t))dy\\
0 & e^{\int_{0}^{\sigma}\nabla\cdot\chi_{\mu}^{\pm}(x(\vartheta,t),t)d\vartheta}-1
\end{bmatrix}\begin{bmatrix}z^{\parallel}(0)\\
z^{\perp}(0)
\end{bmatrix}=-\begin{bmatrix}\Pi^{\parallel}(\sigma)\\
\Pi^{\perp}(\sigma)
\end{bmatrix}-\begin{bmatrix}1\\
0
\end{bmatrix}\dfrac{\sigma_{1}}{\Delta t},\end{split}
\]
where, $e^{\int_{0}^{\sigma}\nabla\cdot\chi_{\mu}^{\pm}(x(\vartheta,t),t)d\vartheta}=\rho_{1}\rho_{2}=\rho_{2}$,
with $\rho_{1}$ and $\rho_{2}$ denoting the Floquet multipliers
\citep{GuckenheimerHolmes1983} of the $\sigma-$periodic limit cycle
$\gamma$ of the ODE (\ref{etamufield}). Solving this system, we
obtain

\begin{equation}
\begin{cases}
z^{\perp}(0)=\dfrac{\Pi^{\perp}(\sigma)}{1-\rho_{2}}\\
\sigma_{1}=\Delta t\Big(\Pi^{\parallel}(\sigma)+z^{\perp}(0)\int_{0}^{\sigma}e^{\int_{0}^{y}\nabla\cdot\chi_{\mu}^{\pm}(x(\vartheta,t),t)d\vartheta}\kappa(x(y,t))dy\Big),
\end{cases}\label{z0andsigma1}
\end{equation}
where the first equation is the same as equation (\ref{z0final}).
The equations in (\ref{z0andsigma1}) are independent of the value
of $z^{\parallel}(0)$ due to the invariance of material flux under
a shift of the parameter $s$. The hyperbolic nature of limit cycles
ensures that $\rho_{1}\rho_{2}=\rho_{2}\neq1$, and thus, both expressions
in (\ref{z0andsigma1}) are well-defined on elliptic OECSs. Observe
that the denominator $(1-\rho_{2})$ is equal to the slope of the
Poincar\'{e} return map along $\gamma$, as shown in \citet{Perko1990}.
The first equation of (\ref{z0andsigma1}) is the only component of
$z(0)$ needed for the computation of the instantaneous material flux
$\varphi_{\gamma(t)}$.

Although $\sigma_{1}$ is not strictly necessary for computing $\varphi_{\gamma(t)}$,
it gives the $\mathcal{O}(\epsilon)$ variation of the period $\sigma_{\epsilon}$
of to the deformed elliptic OECSs, as the parameter $t$ is perturbed
to $t_{\epsilon}=t+\epsilon\Delta t$.

\section{Objectivity of the persistence metric}

\label{Objectivity} Here we show that the non-dimensional metric
$\Theta_{\gamma}(t)$ is objective i.e., invariant under all coordinate
changes of the form 
\begin{equation}
x=Q(t)\tilde{x}+b(t),\label{CoordchangeObj}
\end{equation}
where $Q(t)\in SO(2)$ and $b(t)\in\mathbb{R}^{2}$ are smooth functions
of time. Since the $\Theta_{\gamma}(t)$ is a scalar quantity, in
order for it to be objective \citet{TruesdellNoll2004}, at every point
it must have the same value independent of the actual coordinates
chosen, $x$ or $\tilde{x}$, as long as they are linked by equation
(\ref{CoordchangeObj}). To see this, we check objectivity separately
for the numerator and denominator of (\ref{NondimMetric}).

The spin tensor $W$ introduced in (\ref{VelGradDec}) is well known
to be non-objective \citet{TruesdellNoll2004}, as it transforms as
\[
\tilde{W}=Q^{\top}WQ-Q^{\top}\dot{Q}.
\]
Correspondingly, the plane-normal component $\omega$ of the vorticity
transforms under (\ref{CoordchangeObj}) as 
\[
\tilde{\omega}=\omega-\omega_{Q},
\]
where, $\omega_{Q}$ is such that $Q^{\top}\dot{Q}=\omega_{Q}R$.
The deviation of the vorticity from its spatial mean transforms as
\begin{equation}
\begin{aligned}\tilde{\omega}-\overline{\tilde{\omega}}= & \omega-\omega_{Q}-\dfrac{1}{A_{\partial U}}\int_{A_{\partial U}}(\omega-\omega_{Q})dA\\
= & \omega-\dfrac{1}{A_{\partial U}}\int_{A_{\partial U}}\omega dA=\omega-\overline{\omega},
\end{aligned}
\label{vortdevTrans}
\end{equation}
where, in the second line we used the fact that the domain $U$ is
time independent and $\omega_{Q}$ is space independent. Formula (\ref{vortdevTrans})
proves the objectivity of $\omega_{\gamma(t)}$ defined in (\ref{ObjMeanVort}).

To show the objectivity of the relative material leakage\textit{ } defined
in (\ref{AbsoluteFlux}), we rewrite the pointwise material flux density
(\ref{PointwiseInstFlux}) in the simplified form: 
\begin{equation}
\varphi=\langle\dot{x}_{a1}-\dot{x}_{a2},\ {\Delta x}\rangle,\label{material fluxsimplifiedvObj}
\end{equation}
where, $\dot{x}_{a1}$ and $\dot{x}_{a2}$ represent two general velocity
vectors which have the same base point ${x_{a}}$, and ${\Delta x=x_{a}-x_{b}}$
is a simple distance vector between two points. Representing these
quantities in the $\tilde{x}$ frame, we obtain 
\[
\begin{aligned}\tilde{\dot{x}}_{ai}= & Q^{\top}\dot{x}_{ai}-Q^{\top}\dot{Q}\tilde{x}_{ai}-Q^{\top}\dot{b},\ \ i=1,2,\\
\widetilde{\Delta x}= & Q^{\top}{\Delta x},
\end{aligned}
\]
that, together with (\ref{material fluxsimplifiedvObj}) leads to
\[
\begin{aligned}\tilde{\varphi}= & \langle-Q^{\top}\dot{Q}(\tilde{x}_{a1}-\tilde{x}_{a2})+Q^{\top}(\dot{x}_{a1}-\dot{x}_{a2}),\ Q^{\top}{\Delta x}\rangle,\\
= & \langle Q^{\top}(\dot{x}_{a1}-\dot{x}_{a2}),\ Q^{\top}{\Delta x}\rangle,\\
= & \varphi,
\end{aligned}
\]
where, we used the properties of $Q$ and that $\tilde{x}_{a1}=\tilde{x}_{a2}=\tilde{x}_{a}$
in any coordinate frame since they represent the same base point for
the two velocity vectors involved in the material flux computation.
We have therefore shown that $\omega_{\gamma}$ and $\Gamma_{\gamma}$
are both objective quantities, and hence so is the vortex persistence
metric, $\Theta_{\gamma}$, introduced in Definition \ref{DefVortexPersistence}.

\section{Numerical steps for the computation of elliptic OECSs and $\Theta_{\gamma}$}

\label{SecAutomScheme} 
\begin{algorithm}[H]
\caption{Compute elliptic OECSs \citep{SerraHaller2015}}
\label{algorithm1} \textbf{Input:} A 2-dimensional velocity field. 
\begin{enumerate}
\item Compute the rate-of-strain tensor $S(x,t)=\frac{1}{2}\left(\nabla v(x,t)+\left[\nabla v(x,t)\right]^{T}\right)$
at the current time $t$ on a rectangular grid over the $(x_{1},x_{2}$)
coordinates. 
\item Detect the singularities of $S$ as common, transverse zeros of $S_{11}(\,\cdot\,,t)-S_{22}(\,\cdot\,,t)$
and $S_{12}(\,\cdot\,,t)$, with $S_{ij}$ denoting the entry of $S$
at row $i$ and column $j$. 
\item Determine the type of the singularity (trisector or wedge) as described
in \citet{ShearlessBarrierFarazmand2014}. 
\item Locate isolated wedge-type pairs of singularities and place the Poincar\'{e} 
sections at their midpoint. 
\item Compute the eigenvalue fields $s_{1}(x,t)<s_{2}(x,t)$  and the associated
unit eigenvector fields $e_{i}(x,t)$ of $S(x,t)$ for $i=1,2.$
\item Compute the vector field $\chi_{\mu}^{\pm}(r(s))=\sqrt{\dfrac{s_{2}-\mu}{s_{2}-s_{1}}}e_{1}\ \pm\ \sqrt{\dfrac{\mu-s_{1}}{s_{2}-s_{1}}}e_{2}$
for different values of stretching rate $\mu$, remaining in the range
$\mu\approx0$. 
\item Use the Poincar\'{e}  sections as sets of initial conditions in the computation
of limit cycles of 
\[
x'(s)=\mathrm{sign}\left\langle \chi_{\mu}^{\pm}(x(s)),\tfrac{dx(s-\Delta)}{ds}\right\rangle \chi_{\mu}^{\pm}(x(s)),
\]
where the factor multiplying $\chi_{\mu}^{\pm}(x(s),t)$ removes potential
orientation discontinuities in the direction field $\chi_{\mu}^{\pm}(x(s),t)$
away from singularities, and $\Delta$ denotes the integration step
in the independent variable $s$. 
\end{enumerate}
\textbf{Output:} Elliptic OECSs, related $\chi_{\mu}^{\pm}$ tangent
field and rate of strain tensor field ($S(x,t)$). 
\end{algorithm}

\begin{algorithm}[H]
\caption{Compute the persistence metric for each elliptic OECS}
\label{algorithm2} \textbf{Input:} A 2-dimensional velocity field,
elliptic OECSs, related $\chi_{\mu}^{\pm}$ tangent fields and $S(x,t)$. 
\begin{enumerate}
\item For each elliptic OECS $\gamma$, compute the rotational coherence
$\omega_{\gamma}(t)$. 

\begin{enumerate}
\item[a)] Compute vorticity scalar field $\omega(x,t)$. 
\item[b)] Compute $\omega_{\gamma}(t)$ as: 
\[
\omega_{\gamma}(t)=\dfrac{\vert\int_{A_{\gamma(t)}}[\omega(x,t)-\overline{\omega}(t)]dA\vert}{A_{\gamma(t)}},
\]
where $A_{\gamma(t)}$ is the area enclosed by $\gamma(t)$ and $\overline{\omega}(t)=\tfrac{\int_{A_{\partial U}}\omega(x,t)dA}{A_{\partial U}}$. 
\end{enumerate}
\item For each elliptic OECS, $\gamma$, compute the relative material leakage
, $\Gamma_{\gamma}(t)$: 

\begin{enumerate}
\item[a)] Compute the curvature scalar $\kappa$ and the divergence $\nabla\cdot\chi$
of the $\chi_{\mu}^{\pm}$ tangent field along elliptic OECSs. 
\item[b)] Compute $\partial_{t}S$ along elliptic OECSs using a backward finite
differencing scheme. 
\item[c)] Using equations (\ref{PointwiseInstFlux}-\ref{z0final}), compute
$\Gamma_{\gamma}(t)$ as 
\[
\Gamma_{\gamma}(t)=\dfrac{\oint_{\gamma(t)}\lvert\varphi(x(s,t),t)\rvert ds}{A_{\gamma(t)}}.
\]

\end{enumerate}
\item For each elliptic OECS, $\gamma$, compute $\Theta_{\gamma}(t)=\tfrac{\omega_{\gamma}(t)}{\Gamma_{\gamma}(t)}$. 
\item Within each elliptic OECSs belt (candidate eddy region), select the
one with the maximal $\Theta_{\gamma}(t)$.  
\end{enumerate}
\textbf{Output:} List of coexisting elliptic OECSs with their correspondent
metric value $\Theta_{\gamma}(t)$. 
\end{algorithm}

Here we propose a systematic way to monitor the accuracy of numerical
differentiation involved in equations (\ref{fundmatrsolinetabasesimplified}-\ref{Partialt0etavFt}).
Specifically, equation (\ref{fundmatrsolinetabasesimplified}) requires
spatial differentiation for the computation of $\nabla\cdot\chi$
(step (ii)a), while equation (\ref{Partialt0etavFt}) requires differentiation
in time to compute $\partial_{t}S$ (step (ii)b).

To select the appropriate stepsize for the spatial differentiation
of the $\chi$ field, we turn the relation $(\nabla\chi)^{\top}\chi=0$,
shown in \ref{AppFlux}, into the scalar equation $\langle\chi,\nabla\chi(\nabla\chi)^{\top}\chi\rangle=0$.
The deviation of $\langle\chi,\nabla\chi(\nabla\chi)^{\top}\chi\rangle$
from zero allows to quantitatively monitor the entity of the error
due to spatial differentiation in the material flux computation. For
instance, a complex geometry of the elliptic OECS would require a
finer grid for the accurate computation of $\nabla\chi$. This refinement,
however, is needed only to handle sharp changes in the elliptic OECSs,
which are signaled by high values of the curvature scalar $\kappa$.
Therefore, it is possible to fix a desired maximum allowable deviation
of $\langle\chi,\nabla\chi(\nabla\chi)^{\top}\chi\rangle$ from zero
and select the spatial resolution accordingly.

In a similar fashion, we monitor also the accuracy of numerical finite
differencing in the time direction used to compute $\partial_{t}S$
in (\ref{Partialt0etavFt}). Since the direction field $\chi$ is
normalized, differentiating the identity $\langle\chi_{\mu},\chi_{\mu}\rangle=1$
with respect to time leads to $\partial_{t}\chi_{\mu}\perp\chi_{\mu}$.
Monitoring the deviation of $\vert\langle\tfrac{\partial_{t}\chi_{\mu}}{\lvert\partial_{t}\chi_{\mu}\rvert},\chi_{\mu}\rangle\vert$
from zero allows a systematic assessment of the appropriate time step
required to compute $\partial_{t}S$.

The quantities $\langle\chi,\nabla\chi(\nabla\chi)^{\top}\chi\rangle$
and $\vert\langle\tfrac{\partial_{t}\chi_{\mu}}{\lvert\partial_{t}\chi_{\mu}\rvert},\chi_{\mu}\rangle\vert$
play the role of numerical reliability parameters and allow us to
compute the material flux through any elliptic OECS in an efficient
and accurate fashion.

\section{Ocean surface flow dataset}

\label{SSHdataset} Under the geostrophic assumption, the ocean surface
topology measured by satellites plays the role of a stream function
for the related surface currents. With $h$ denoting the sea surface
height, the velocity field in longitude-latitude coordinates $[\phi,\ \theta]$,
can be expressed as 
\[
\dot{\phi}=-\dfrac{g}{R^{2}f(\theta)\cos\theta}\partial_{\theta}h(\phi,\theta,t),\ \ \ \text{and}\ \dot{\theta}=\dfrac{g}{R^{2}f(\theta)\cos\theta}\partial_{\phi}h(\phi,\theta,t),
\]
where $f(\theta):=2\Omega\sin\theta$ denotes the Coriolis parameter,
$g$ the constant of gravity, $R$ the mean radius of the earth and
$\Omega$ its mean angular velocity. The velocity field is available
at weekly intervals, with a spatial longitude-latitude resolution
of $0.25^{\circ}$. For more detail on the data, see \citet{Beron-VeraDatassh2013}.

\section{Thresholding requirement for common Eulerian diagnostics}

\label{OWThresholding}
\begin{figure}
\centering \includegraphics[width=0.9\textwidth]{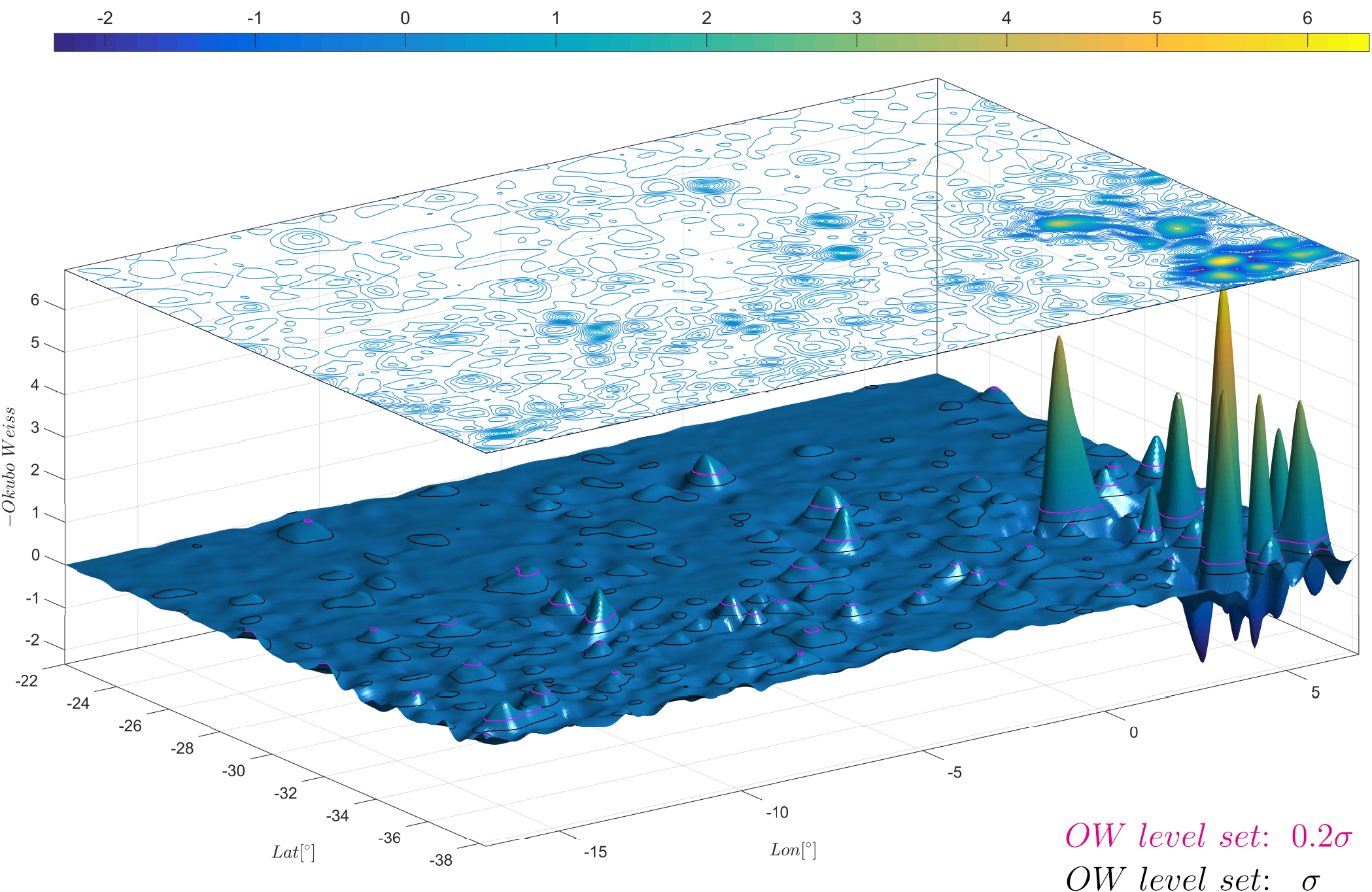}
\caption{$OW$ parameter and two specified level sets corresponding to $OW=0.2\sigma$
(magenta) and $OW=\sigma$ (black) with $\sigma$ being the $OW$
spatial standard deviation.}
\label{OWans2levsets} 
\end{figure}

Vortex definitions based on scalar fields (e.g.,
$OW$) are often ambiguous due to their dependence on ad hoc of thresholding
parameters. For the $OW$-criterion, this threshold value is typically
$\alpha\sigma$, with $\sigma$ being the spatial standard deviation
of the $OW$ parameter, and $\alpha\in\mathbb{R}$ selected as a problem-dependent
constant. Figure \ref{OWans2levsets} shows the $OW$ level sets for
two different values of $\alpha$: 0.2 and 1, as suggested in \citet{OWls_theshld_Henson2008}
and \citet{no0lsOW_Koszalka2009}, respectively. Note how the values
of $\alpha$ can significantly change the overall number and geometry
of vortices identified.

\begin{figure}[t]
\centering \includegraphics[width=0.95\textwidth]{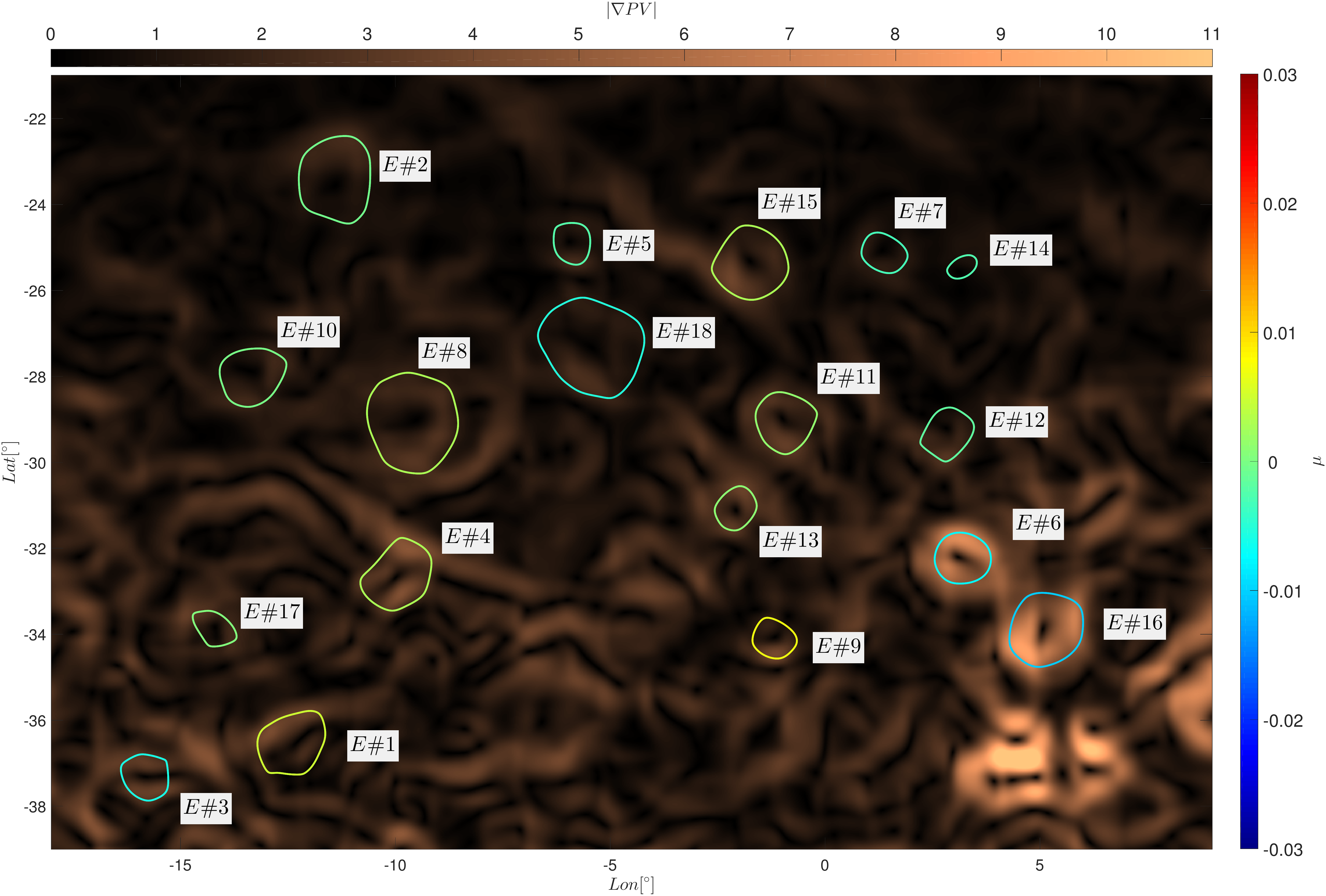}
\caption{Elliptic OECSs with highest vortex persistence metric $\Theta_{\gamma}$,
plotted over $\vert\nabla PV\vert$ (horizontal colorbar). Elliptic
OECSs are encoded with a color representing their stretching-rate
value $\mu$ (right colorbar). Black numbers label different vortical
regions encircled by elliptic OECSs.}
\label{BestISGonPVGrad}
\end{figure}

In Fig. \ref{BestISGonPVGrad}, we show the elliptic OECSs with highest
vortex persistence metric $\Theta_{\gamma}$, shown in Fig. \ref{BestISGonOW3dsurf},
on a scalar field representing the $\vert\nabla PV\vert$ where $PV$
is computed as in \citet{Early2011}. Regions of high $PV$ gradient
are frequently used indicators of instantaneous ellipticity in unsteady
fluid flows. In the south-east of the domain ($[1{}^{\circ}E,7^{\circ}E]$,
$[31^{\circ}S,38^{\circ}S]$), although there are several regions
of high $\vert\nabla PV\vert$, only one long-lived Lagrangian eddy
is present. At the same time, $\vert\nabla PV\vert$ fails to signal
several other regions captured by elliptic OECSs (see e.g., eddies
$\#2,\ \#8,\text{ {and} }\#18$), where long-lived Lagrangian eddies
are present. Therefore, a prediction based only on $\vert\nabla PV\vert$,
i.e., choosing an ad hoc threshold parameter to locate vortices instead
of using elliptic OECSs, would be even weaker than the one shown in
Fig. \ref{MetricGLMOW}c. A similar conclusion holds for the $OW$-criterion,
as discussed in section \ref{PredonOcean}. 
\newpage{} 

\bibliographystyle{plainnat}
\bibliography{ReferenceList2_Local}

\begin{thebibliography}{27}
\providecommand{\natexlab}[1]{#1}
\providecommand{\url}[1]{\texttt{#1}}
\expandafter\ifx\csname urlstyle\endcsname\relax
  \providecommand{\doi}[1]{doi: #1}\else
  \providecommand{\doi}{doi: \begingroup \urlstyle{rm}\Url}\fi

\bibitem[Arnold(1973)]{ArnoldODE1973}
V.~Arnold.
\newblock \emph{Ordinary Differential Equations}.
\newblock MIT Press, Boston, 1973.

\bibitem[Batchelor(2000)]{Batchelor2000}
G.~K. Batchelor.
\newblock \emph{{An introduction to fluid dynamics}}.
\newblock Cambridge univ. press, 2000.

\bibitem[Beal et~al.(2011)Beal, De~Ruijter, Biastoch, Zahn, and
  GROUP]{BealNature2011}
L.~M. Beal, W.~P.~M. De~Ruijter, A.~Biastoch, R.~Zahn, and
  SCOR/WCRP/IAPSO~WORKING GROUP.
\newblock {On the role of the Agulhas system in ocean circulation and climate}.
\newblock \emph{Nature}, 472\penalty0 (7344):\penalty0 429--436, 2011.

\bibitem[Beron-Vera et~al.(2013)Beron-Vera, Wang, Olascoaga, Goni, and
  Haller]{Beron-VeraDatassh2013}
F.~J. Beron-Vera, Y.~Wang, M.~J. Olascoaga, G.~J. Goni, and G.~Haller.
\newblock {Objective detection of oceanic eddies and the Agulhas leakage}.
\newblock \emph{J. Phys. Oceanogr.}, 43:\penalty0 1426--1438, 2013.

\bibitem[Chelton et~al.(2007)Chelton, Schlax, Samelson, and
  De~Szoeke]{Chelton2007}
D.~B. Chelton, M.~G. Schlax, R.~M. Samelson, and R.~A. De~Szoeke.
\newblock {Global observations of large oceanic eddies}.
\newblock \emph{Geophys. Res. Lett.}, 34, 2007.

\bibitem[Early et~al.(2011)Early, Samelson, and Chelton]{Early2011}
J.~J. Early, R.~M. Samelson, and D.~B. Chelton.
\newblock {The evolution and propagation of quasigeostrophic ocean Eddies*}.
\newblock \emph{J. Phys. Oceanogr.}, 41:\penalty0 1535--1555, 2011.
\newblock \doi{10.1175/2011JPO4601.1}.

\bibitem[Farazmand and Haller(2016)]{Farazmand2015}
M.~Farazmand and G.~Haller.
\newblock {Polar rotation angle identifies elliptic islands in unsteady
  dynamical systems}.
\newblock \emph{Physica D}, 315:\penalty0 1 -- 12, 2016.
\newblock ISSN 0167-2789.
\newblock \doi{http://dx.doi.org/10.1016/j.physd.2015.09.007}.
\newblock URL
  \url{http://www.sciencedirect.com/science/article/pii/S0167278915001748}.

\bibitem[Farazmand et~al.(2014)Farazmand, Blazevski, and
  Haller]{ShearlessBarrierFarazmand2014}
M.~Farazmand, D.~Blazevski, and G.~Haller.
\newblock {Shearless transport barriers in unsteady two-dimensional flows and
  maps}.
\newblock \emph{Physica D}, 278:\penalty0 44--57, 2014.

\bibitem[Griffa et~al.(2007)Griffa, Kirwan, Mariano, {\"O}zg{\"o}kmen, and
  Rossby]{Griffa2007}
A.~Griffa, A.~D. Kirwan, A.~J. Mariano, T.~{\"O}zg{\"o}kmen, and H.~T. Rossby.
\newblock \emph{{Lagrangian analysis and prediction of coastal and ocean
  dynamics}}.
\newblock Cambridge University Press, 2007.

\bibitem[Guckenheimer and Holmes(1983)]{GuckenheimerHolmes1983}
J.~Guckenheimer and P.~Holmes.
\newblock \emph{Nonlinear oscillations, dynamical systems, and bifurcations of
  vector fields}, volume~42.
\newblock Springer Science \& Business Media, 1983.

\bibitem[Haller(2015)]{LCSHallerAnnRev2015}
G.~Haller.
\newblock {Lagrangian coherent structures}.
\newblock \emph{Annual Rev. Fluid. Mech}, 47:\penalty0 137--162, 2015.

\bibitem[Haller and Beron-Vera(2013)]{BlackHoleHaller2013}
G.~Haller and F.~J. Beron-Vera.
\newblock {Coherent Lagrangian vortices: the black holes of turbulence}.
\newblock \emph{J. Fluid Mech.}, 731, 9 2013.
\newblock ISSN 1469-7645.
\newblock \doi{10.1017/jfm.2013.391}.
\newblock URL \url{http://journals.cambridge.org/article_S0022112013003911}.

\bibitem[Haller and Iacono(2003)]{HallerIacono2003}
G.~Haller and R.~Iacono.
\newblock {Stretching, alignment, and shear in slowly varying velocity fields}.
\newblock \emph{Phys. Rev. E}, 68:\penalty0 056304, 2003.

\bibitem[Haller et~al.(2016)Haller, Hadjighasem, Farazmand, and
  Huhn]{HallerLAVD2015}
G.~Haller, A.~Hadjighasem, M.~Farazmand, and F.~Huhn.
\newblock {Defining coherent vortices objectively from the vorticity}.
\newblock \emph{J. Fluid Mech.}, 795:\penalty0 136--173, 5 2016.
\newblock ISSN 1469-7645.
\newblock \doi{10.1017/jfm.2016.151}.
\newblock URL \url{http://journals.cambridge.org/article_S0022112016001518}.

\bibitem[Henson and Thomas(2008)]{OWls_theshld_Henson2008}
S.~A. Henson and A.~C. Thomas.
\newblock {A census of oceanic anticyclonic eddies in the Gulf of Alaska}.
\newblock \emph{Deep Sea Res.Part I: Oceanogr. Res. Papers}, 55:\penalty0
  163--176, 2008.

\bibitem[Isern-Fontanet et~al.(2003)Isern-Fontanet, Garc{\'\i}a-Ladona, and
  Font]{Isern-Fontanet2003}
J.~Isern-Fontanet, E.~Garc{\'\i}a-Ladona, and J.~Font.
\newblock {Identification of marine eddies from altimetric maps}.
\newblock \emph{J. of Atmosph. and Oceanic Technology}, 20:\penalty0 772--778,
  2003.

\bibitem[Isern-Fontanet et~al.(2004)Isern-Fontanet, Font, Garc{\'\i}a-Ladona,
  Emelianov, Millot, and Taupier-Letage]{Isern-Fontanet2004}
J.~Isern-Fontanet, J.~Font, E.~Garc{\'\i}a-Ladona, M.~Emelianov, C.~Millot, and
  I.~Taupier-Letage.
\newblock {Spatial structure of anticyclonic eddies in the Algerian basin
  (Mediterranean Sea) analyzed using the Okubo--Weiss parameter}.
\newblock \emph{Deep Sea Res. Part II: Topical Studies in Oceanogr.},
  51:\penalty0 3009--3028, 2004.

\bibitem[Isern-Fontanet et~al.(2006)Isern-Fontanet, Garc{\'\i}a-Ladona, and
  Font]{OWlsthrsh_Isern-Fontanet2006}
J.~Isern-Fontanet, E.~Garc{\'\i}a-Ladona, and J.~Font.
\newblock {Vortices of the Mediterranean Sea: An altimetric perspective}.
\newblock \emph{J. Phys. Oceanogr.}, 36:\penalty0 87--103, 2006.

\bibitem[Karrasch et~al.(2015)Karrasch, Huhn, and Haller]{AutoMdetectFlorian}
D.~Karrasch, F.~Huhn, and G.~Haller.
\newblock {Automated detection of coherent Lagrangian vortices in
  two-dimensional unsteady flows}.
\newblock In \emph{Proc. R. Soc. Lond. A.}, volume 471. The Royal Society,
  2015.

\bibitem[Koszalka et~al.(2009)Koszalka, Bracco, McWilliams, and
  Provenzale]{no0lsOW_Koszalka2009}
I.~Koszalka, A.~Bracco, J.~C. McWilliams, and A.~Provenzale.
\newblock {Dynamics of wind-forced coherent anticyclones in the open ocean}.
\newblock \emph{J. Geophys. Res.: Oceans (1978--2012)}, 114, 2009.

\bibitem[Okubo(1970)]{Okubo1970}
A.~Okubo.
\newblock {Horizontal dispersion of floatable particles in the vicinity of
  velocity singularities such as convergences}.
\newblock In \emph{Deep-Sea Res.}, volume~17, pages 445--454. Elsevier, 1970.

\bibitem[Perko(1990)]{Perko1990}
L.~M. Perko.
\newblock {Global families of limit cycles of planar analytic systems}.
\newblock \emph{Trans. of the American Math. Society}, 322:\penalty0 627--656,
  1990.

\bibitem[Serra and Haller(2016)]{SerraHaller2015}
M.~Serra and G.~Haller.
\newblock {Objective Eulerian Coherent Structures}.
\newblock \emph{Chaos, (in press)}, 2016.

\bibitem[Truesdell and Noll(2004)]{TruesdellNoll2004}
C.~Truesdell and W.~Noll.
\newblock \emph{{The non-linear field theories of mechanics}}.
\newblock Springer, 2004.

\bibitem[Wang et~al.(2016)Wang, Beron-Vera, and Olascoaga]{Wang2016}
Y~Wang, F.~J. Beron-Vera, and M.~J. Olascoaga.
\newblock {The life cycle of a coherent Lagrangian Agulhas ring}.
\newblock \emph{arXiv preprint arXiv:1601.01560}, 2016.

\bibitem[Waugh et~al.(2006)Waugh, Abraham, and Bowen]{Waugh2006}
D.~W. Waugh, E.~R. Abraham, and M.~M. Bowen.
\newblock {Spatial variations of stirring in the surface ocean: A case study of
  the Tasman Sea}.
\newblock \emph{J. Phys. Oceanogr.}, 36:\penalty0 526--542, 2006.

\bibitem[Weiss(1991)]{Weiss1991}
J.~Weiss.
\newblock {The dynamics of enstrophy transfer in two-dimensional
  hydrodynamics}.
\newblock \emph{Physica D}, 48:\penalty0 273--294, 1991.

\end{thebibliography}

\end{document}